\documentclass[12pt]{article}
\usepackage{amsmath,amssymb,amsthm, amscd, pb-diagram}

\begin{document}

\begin{center}
{\bf Towards a characterization of toric hyperk\"{a}hler varieties among symplectic singularities II}
\end{center}
\vspace{0.4cm}

\begin{center}
{\bf Yoshinori Namikawa}      
\end{center}
\vspace{0.2cm}

\begin{center}
{\bf \S 1. Introduction}. 
\end{center}

This is a continuation of [Part I] with the same title. We use the same notation and convention as in [Part I]. 
In [Part I] we have proved the following.  \vspace{0.2cm} 

{\bf Theorem 5.8}. {\em Let $(X, \omega)$ be a conical symplectic variety of dimension $2n$ which has a projective symplectic resolution. 
Assume that $X$ admits an effective Hamiltonian action of an $n$-dimensional algebraic torus $T^n$, compatible with the conical 
$\mathbf{C}^*$-action. Then there is a $T^n$-equivariant (complex analytic) isomorphism 
$\varphi: (X, \omega) \to (Y(A, 0), \omega_{Y(A,0)})$ 
which makes the following diagram commutative}  
\begin{equation} 
\begin{CD} 
(X, \omega) @>{\varphi}>> (Y(A, 0), \omega_{Y(A,0)}) \\  
@V{\mu}VV @V{\bar{\mu}}VV  \\  
(\mathfrak{t}^n)^* @>{id}>> (\mathfrak{t}^n)^*.
\end{CD} 
\end{equation}   
{\em Here $A$ is unimodular and the vertical maps are moment maps for the $T^n$-actions.} 
{\em Moreover,  we have $\varphi (0_X) = 0_{Y(A, 0)}$.} 
\vspace{0.2cm}

In this article we prove that, when $wt(\omega) = 2$, after replacing the original conical $\mathbf{C}^*$-action on $(X, \omega)$ by another  (algebraic) conical $\mathbf{C}^*$-action on $(X, \omega)$, we can take a new isomorphism $\varphi'$ in a $\mathbf{C}^*$-equivariant way. In particular, such $\varphi'$ is an algebraic isomorphism. This answers Question 5.10 posed in [Part I].

The following is a basic strategy for the proof. In the situation of Theorem 5.8 of [Part I], we pull back the standard conical $\mathbf{C}^*$-action on 
$(Y(A, 0), \omega_{Y(A,0)})$ by $\varphi$ and get a conical $\mathbf{C}^*$-action on $(X, \omega)$. However, this is a complex analytic $\mathbf{C}^*$-action. Here we replace the original $\mathbf{C}^*$-action on $X$ with a new algebraic $\mathbf{C}^*$-action on $X$ by  using  the Hamiltonian $T^n$-action (cf. (2.2)) and compare it with the complex analytic $\mathbf{C}^*$-action.  
Then we can construct a complex analytic $T^n$-equivariant automorphism $\phi: (X, \omega) \to (X, \omega)$ so that the 
complex analytic $\mathbf{C}^*$-action on the lefthand side is transformed to the new $\mathbf{C}^*$-action on the righthand side. Now the composite $\varphi' := \varphi \circ \phi^{-1}$ is the desired one.

The precise statement of our main theorem is:  
\vspace{0.2cm}

{\bf Theorem}. {\em Let $(X, \omega)$ be a conical symplectic variety of dimension $2n$ with $wt(\omega) = 2$, which has a projective symplectic resolution. Assume that $X$ admits an effective Hamiltonian action of an $n$-dimensional algebraic torus $T^n$, compatible with the conical $\mathbf{C}^*$-action. Replace the conical $\mathbf{C}^*$-action by a suitable new conical $\mathbf{C}^*$-action compatible with the $T^n$-action.   
Then one can find a $T^n$-equivariant isomorphism 
$\varphi: (X, \omega) \to (Y(A, 0), \omega_{Y(A,0)})$ 
so that the following diagram commutes   
\begin{equation} 
\begin{CD} 
(X, \omega) @>{\varphi}>> (Y(A, 0), \omega_{Y(A,0)}) \\  
@V{\mu}VV @V{\bar{\mu}}VV  \\  
(\mathfrak{t}^n)^* @>{id}>> (\mathfrak{t}^n)^*.
\end{CD} 
\end{equation}  
and $\varphi$  is $\mathbf{C}^*$-equivariant. 
Here $A$ is unimodular and the vertical maps are moment maps for the $T^n$-actions.}
\vspace{0.2cm}

The theorem has been also announced by A. Hubbard in [Hu]. His proof uses Cox rings and is quite different from ours. 
The author came up with the present proof in March, 2026.
\vspace{0.2cm}
   
\begin{center}
{\bf \S 2. Proof of Theorem}. 
\end{center}

We start with the situation of Theroem 5.8 of [Part I].

(2.1)   Let $\pi: \tilde{X} \to X$ be a projective symplectic resolution. 
The Hamiltonian $T^n$-action on $(X, \omega)$ lifts to a Hamiltonian $T^n$-action on $(\tilde{X}, \omega_{\tilde X})$. 
Put $\tilde{\mu} := \mu \circ \pi$. We shall construct a $T^n$-invariant Zariski open subset $U$ of $\tilde{X}$ with the
following property:  \vspace{0.2cm}

(*)  $\tilde{\mu}\vert_U: U \to (\mathfrak{t}^n)^*$ is a principal $T^n$-bundle.  
 \vspace{0.2cm}
   
Let us consider the original (algebraic) $\mathbf{C}^*$-action on $X$.  It also lifts to a $\mathbf{C}^*$-action on $\tilde{X}$. 
Since these two actions commute, we have a $\mathbf{C}^* \times T^n$-action on $\tilde{X}$. 
By Sumihiro's theorem, $\tilde{X}$ is covered by finitely many $\mathbf{C}^* \times T^n$-invariant affine subset $U_i$: 
$$ \tilde{X} = \cup_{i \in I} U_i.$$ 
Restrict the map $\tilde{\mu}$ to $U_i$: $$\mu_i: U_i \to (\mathfrak{t}^n)^*.$$ Then $\mu_i$ is the moment map for the $T^n$-action on 
$(U_i, \omega_{\tilde X}\vert_{U_i})$. By [Part I, Theorem 1.3, (1)], $\mu_i$ factorizes as 
$$U_i \stackrel{\tau_i}\to U_i/\hspace{-0.1cm}/T^n \stackrel{\nu_i}\to (\mathfrak{t}^n)^*$$ and $\nu_i$ is an etale map. 
Since $U_i$ is $\mathbf{C}^*$-invariant, we see that $\nu_i$ is an isomorphism and $\mu_i$ is a $\mathbf{C}^*$-equivariant map by  
[Part I, Corollary 3.2]. Here the $\mathbf{C}^*$-action on $(\mathfrak{t}^n)^*$ is given by the scaling action 
$\times \sigma^2$ with $\sigma \in \mathbf{C}^*$ because $wt(\omega) = 2$ (cf. [Part I, Proposition 3.1]). 

By [Lo, Lemma 2.6] one can take a $T^n$-invariant Zariski open subset $U_i^0$ of $U_i$ in such a way that 
$\mu_i\vert_{U_{i^0}} : U_i^0 \to (\mathfrak{t}^n)^*$ is a principal $T^n$-bundle. 

By the local description of $\mu_i$ (cf. [Part I, Example 1.2]), there are only finitely many choices of such $U_i^0$. 
Since $\tilde{X}$ is covered by finitely many $U_i$, there are only finitey many $T^n$-invariant Zariski open subsets 
of $\tilde{X}$ with the property (*). 

Let $U$ be one of such $(U_i)^0$. Then $U$ must be preserved by the $\mathbf{C}^*$-action. 
In fact, $\sigma \cdot U$ ($\sigma \in \mathbf{C}^*$) is also a $T^n$-invariant open set of $\tilde{X}$ with the property (*). 
But, as remarked above, there are only finitely many such open sets. Since $\mathbf{C}^*$ is connected, we have 
$\sigma \cdot U = U$.   
\vspace{0.2cm}

(2.2) Given $\varphi$ in Theorem 5.8 of [Part I], we consider the 2-nd $\mathbf{C}^*$-action on $X$ induced from the standard $\mathbf{C}^*$-action on $Y(A,0)$ by $\varphi$. This $\mathbf{C}^*$-action also lifts to a $\mathbf{C}^*$-action on $\tilde{X}$. To distinguish this 
$\mathbf{C}^*$-action from the orginal $\mathbf{C}^*$-action on $\tilde{X}$, we denote the original $\mathbf{C}^*$-action by 
$$ \cdot: \mathbf{C}^* \times \tilde{X} \to \tilde{X}, \:\:\: (\sigma, z) \to \sigma \cdot z $$
and denote the 2-nd $\mathbf{C}^*$-action by  
$$ \bullet : \mathbf{C}^* \times \tilde{X} \to \tilde{X}, \:\:\: (\sigma, z) \to \sigma \bullet z. $$
Let $U$ be a $T^n$-invariant Zariski open subset of $\tilde{X}$ with  property (*) (cf. (2.1)). Then, since the 2-nd $\mathbf{C}^*$-action commutes  
with the $T^n$-action, we see that  $U$ is preserved by the 2-nd $\mathbf{C}^*$-action by the same reasoning 
as in (2.1).   
As a consequence, $U$ has 2 different $\mathbf{C}^*$-actions, each of them commutes with the $T^n$-action on $U$: 
$$\cdot : \mathbf{C}^* \times U \to U, \:\:\: (\sigma, z) \to \sigma \cdot z, \:\:\:\:\:\: \bullet : \mathbf{C}^* \times U \to U, \:\:\: (\sigma, z) \to \sigma \bullet z$$ Let us simply write $\mu_U$ for $\tilde{\mu}\vert_U$.
Notice that $\mu_U (\sigma \cdot z) = \mu_U(\sigma \bullet z)$ for $(\sigma, z) \in \mathbf{C}^* \times U$ because, for both $\mathbf{C}^*$-actions,  $wt(\omega) = 2$ and $\mu_U$ is $\mathbf{C}^*$-equivariant. 
 We compare these two actions. Since $U$ is a principal $T^n$-bundle over $(\mathfrak{t}^n)^*$, there is a holomorphic map 
 $h: \mathbf{C}^* \times U \to T^n, \:\:\: (\sigma, z) \to h(\sigma, z)$ such that $$\sigma \bullet z = h(\sigma, z)\sigma \cdot z.$$
  
 By using the fact that these two actions respectively commute with the $T^n$-action and $T^n$ is commutative, one sees that 
 $h$ factorizes as $$\mathbf{C}^* \times U \stackrel{id \times \mu_U}\to \mathbf{C}^* \times (\mathfrak{t}^n)^* \stackrel{\bar{h}}\to 
 T^n$$ for a suitable holomorphic map $\bar{h}: \mathbf{C}^* \times (\mathfrak{t}^n)^* \to T^n$. 
 In fact, for $t \in T^n$, we have $$\sigma \bullet tz = t \sigma \bullet z = t h(\sigma, z) \sigma \cdot z = h(\sigma, z)t \sigma\cdot z 
 = h(\sigma, z)\sigma \cdot (tz).$$ Therefore, $h(\sigma, z) = h(\sigma, tz)$ and $h(\sigma, z)$ depends only on $\mu_U(z)$. 
 
Moreover, assume that $\mu_U(z) = 0$. Then, for any $\sigma \in \mathbf{C}^*$, we have $\mu_U(\sigma \cdot z) = \sigma \mu_U(z) = \sigma \cdot 0 = 0$ and $\mu_U(\sigma \bullet z) = 0$. 
Now, for any $\sigma, \tau \in \mathbf{C}^*$, we have 
$$(\tau \sigma)\bullet z = \tau \bullet (\sigma \bullet z) = \bar{h}(\tau, 0) \tau \cdot (\sigma \bullet z)  = $$
$$\bar{h}(\tau, 0)\tau \cdot \bar{h}(\sigma, 0)\sigma \cdot z = \bar{h}(\tau, 0)\bar{h}(\sigma, 0)(\tau \sigma)\cdot z.$$ 
This means that $$\bar{h}(\tau\sigma, 0) = \bar{h}(\tau, 0)\bar{h}(\sigma, 0).$$ 
Therefore, we see that $\iota := \bar{h}\vert_{\mathbf{C}^* \times \{0\}} : \mathbf{C}^* \to T^n$ is a 
homomorphism of complex Lie groups.

Now $\sigma \in \mathbf{C}^*$ acts on $U$ by $$u \to \sigma \cdot \iota (\sigma) \cdot u.$$ 
We denote this new $\mathbf{C}^*$-action by 
$$ * : \mathbf{C}^* \times U \to U, \:\:\: (\sigma, z) \to \sigma * z.$$
Notice that this $\mathbf{C}^*$-action on $U$ extends to a $\mathbf{C}^*$-action on $\tilde{X}$, which descends to a 
$\mathbf{C}^*$-action on $X$. They are all algebraic action by definition. 
By the definition of $*$ we also have 
$$ *\vert_{\mathbf{C}^* \times \mu_U^{-1}(0)} = \bullet\vert_{\mathbf{C}^* \times \mu_U^{-1}(0)}.$$      
 
(2.3) The $\mathbf{C}^*$-actions $*$ and $\bullet$ respectively determine vector fields $v$ and $v'$ on $U$. 
Put $\omega_U := \omega_{\tilde X}\vert_U$. 
We prove \vspace{0.2cm} 

{\bf Proposition 1}. {\em There is a complex analytic $T^n$-equivariant automorphism $\phi$ of $(U, \omega_U)$ such that 
$\phi_*v' = v$ and 
the following diagram commutes:} 
\begin{equation} 
\begin{CD} 
(U, \omega_U) @>{\phi}>> (U, \omega_U) \\  
@V{\mu_U}VV @V{\mu_U}VV  \\  
(\mathfrak{t}^n)^* @>{id}>> (\mathfrak{t}^n)^*.
\end{CD} 
\end{equation}   
 
{\em Proof}. Put $\eta := v \rfloor \omega_U$ and $\eta' := v' \rfloor \omega_U$. It suffices to find such a $\phi$ that 
satisfies $\phi^*\eta = \eta'$. In fact, if $\phi^*\eta = \eta'$, then  
$$\phi^*(v \rfloor \omega_U) =  v' \rfloor \omega_U = v' \rfloor \phi^*\omega_U = \phi^*(\phi_*v' \rfloor \omega_U).$$ 
This means that $v \rfloor \omega_U = \phi_*v' \rfloor \omega_U$, which implies that $v = \phi_*v'$ because $\omega_U$ is 
nondegenerate.  
By assumption we have $L_v \omega_U = 2\omega_U$, where the left hand side equals $d(v \rfloor \omega_U) + v \rfloor d\omega_U$.  
Since $d\omega_U = 0$, we see that $d\eta = 2\omega_U$. Similarly, we have $d\eta' = 2\omega_U$. 
We put $\eta_t :=  t\eta + (1-t)\eta'$ for $t \in [0,1]$. Then 
$$\frac{d\eta_t}{dt} = \eta - \eta.'$$ We note here that $(\mu_U)_*v = (\mu_U)_*v' = 2\zeta$, where $\zeta$ is the Euler 
vector field on $(\mathfrak{t}^n)^*$: 
$$ \zeta = \sum t_i \frac{\partial}{\partial t_i}, \:\:\: (t_1, ..., t_n) \in (\mathfrak{t}^n)^*.$$       
In particular, $v - v' \in \Theta_{U/(\mathfrak{t}^n)^*}$. Since $\mu_U$ is a Lagrangian fibration with respect to $\omega_U$, 
we see that $\eta - \eta' \in (\mu_U)^*\Omega^1_{(\mathfrak{t}^n)^*}$. Since $\eta - \eta'$ is $T^n$-invariant, we can write 
$\eta - \eta' = (\mu_U)^*\beta$ for some $\beta \in \Omega^1_{(\mathfrak{t}^n)^*}$. Since $\eta - \eta'$ is d-closed, the 1-form $\beta$ is a 
d-closed 1-form on $(\mathfrak{t}^n)^*$; hence we can write $\beta = df$ for some holomorphic function $f$ on      
$(\mathfrak{t}^n)^*$. There is one more constraint for $f$. In fact, as we noticed in the last part of (2.2), we have $(v - v')\vert (\mu_U)^{-1}(0)$. This implies that $(\eta - \eta')\vert (\mu_U)^{-1}(0) = 0$. Therefore we have $df(0) = 0$. If we write $f$ as a powerseries 
$f = f_0 + f_1 + f_2 + ... $ with $\mathrm{deg}(f_i) = i$, then we may assume that $f_0 = f_1 = 0$.  
As a consequence, we can write 
$$\eta - \eta' = (\mu_U)^*df \:\: \mathrm{with} \:\: f = f_2 + f_3 + ... $$  
By using Moser's trick, we shall find a family $\{\phi_t\}_{0 \le t \le 1}$ of $T^n$-invariant complex analytic automorphisms of 
$(U, \omega_U)$ such that 
$$\frac{d}{dt}(\phi_t^*\eta_t) = 0.$$ 
For a function $g$ on $(\mathfrak{t}^n)^*$, we put $X := H_{(\mu_U)^*g}$. 
Let $\{\phi_t\}$ be the family of automorphisms $(U, \mu_U)$ determined by $X$. Namely they satisfy 
$$\frac{d\phi_t}{dt} = X(\phi_t), \:\:\: \phi_0 = id_U.$$
Then we have 
$$\frac{d}{dt}(\phi_t^*\eta_t) = \phi_t^*(L_X\eta_t + \frac{d\eta_t}{dt}) = 
\phi_t^*(d(X \rfloor \eta_t) + X \rfloor d\eta_t + \mu_U^*df).$$
Since $d\eta_t = t(d\eta - d\eta') + d\eta = t(2\omega_U -2 \omega_U) + 2\omega_U = 2\omega_U$, 
we have 
$$X \rfloor d\eta_t = X \rfloor 2\omega_U = 2H_{(\mu_U)^*g} \rfloor \omega_U = 2\omega_U(H_{(\mu_U)^*g},  \cdot) 
= -2 d\mu_U^*g.$$
Put $v_t := tv + (1-t)v'$. Then $\eta_t = \omega_U(v_t, \cdot)$.  Then we have 
$$H_{(\mu_U)^*g}   \rfloor \eta_t = \omega_(v_t, H_{(\mu_U)^*g}) = d\mu_U^*g (v_t) = v_t \rfloor \mu_U^*dg 
= \mu_U^*((\mu_U)_*v_t \rfloor dg) = 2\mu_U^*(\zeta \rfloor dg).$$ In the last equality, we use the fact that $\mu_U)_*v_t = 2\zeta$. 
Therefore we  have 
$$\frac{d}{dt}(\phi_t^*\eta_t) = \phi_t^*(\mu_U^*d(2\zeta (dg) - 2g + f)).$$ 
We can solve the equation $$2\zeta (dg) - 2g  = -f$$ for $g$. In fact, let us write 
$g = g_0 + g_1 + g_2 + ... $, then $2\zeta (dg) - 2g  = \sum 2(i -1)g_i$. 
When $i \ne 1$, we just put $g_i := \frac{-1}{2(i-1)}f_i$. When $i = 1$, we can take an arbitrary $g_1$ because $f_1 = 0$ (This part is crucial). 
If we put $\phi := \phi_1$, then $\eta' = \phi^*\eta$. 
$\square$ 
\vspace{0.2cm}

(2.4) In [Part I] we have defined a closed subset $F_X \subset (\mathfrak{t}^n)^*$ with $\mathrm{Codim}_{((\mathfrak{t}^n)^*} F_X 
\geq 2$ and have  described the local structure of the moment map $\mu: X \to (\mathfrak{t}^n)^*$ around 
each  point $t \in (\mathfrak{t}^n)^* - F_X$. Put $(\mathfrak{t}^n)^{*, 0} := (\mathfrak{t}^n)^* - F_X$ and 
$X^0 := \mu^{-1}((\mathfrak{t}^n)^{*, 0})$. 
For a holomorphic function $g$ on $(\mathfrak{t}^n)^{*, 0}$, one can integrate $H_{\mu^*g}$ to  $\mathrm{exp}(H_{\mu^*g})$, which 
is a $T^n$-invariant symplectic automorphism of $(X^0, \omega\vert_{X^0})$ over $(\mathfrak{t}^n)^{*, 0}$. 
(See [Lo], [Part I, Introduction and p. 42] for details.) This means that the automorphism $\phi$ in Proposition 1 determines an automorphism of $X^0$. By Proposition 2,1, (2) of [Part I], every fiber of $\mu: X \to (\mathfrak{t}^n)^*$ has dimension  $n$. Therefore, 
we have $\mathrm{Codim}_{X}(X - X^0) \geq 2$. Since $X$ is a Stein normal variety, this automorphism of $X^0$ uniquely extends to an 
automorphism of $X$, which we denote by the same $\phi$. By the construction, $\phi : (X, \omega) \to (X, \omega)$ transforms 
the $\mathbf{C}^*$-action $\bullet$ on the lefthand side to the $\mathbf{C}^*$-action $*$ on the righthand side. Then 
$$\varphi \circ \phi^{-1}   : (X, \omega) \stackrel{\phi^{-1}}\to (X, \omega) \stackrel{\varphi}\to (Y(A, 0), \omega_{Y(A,0)})$$ 
transforms the $\mathbf{C}^*$-action $*$ to the standard $\mathbf{C}^*$-action on $Y(A, 0)$. Note that $\phi^{-1}$ and 
$\varphi$ are only complex analytic isomorphisms, but the composite $\varphi \circ \phi^{-1}$ is an algebraic isomorphism. 
This completes the proof of Theorem.  
\vspace{0.5cm}

\begin{center}
{\bf References}
\end{center}
 
[Part I] Namikawa, Y.: Towards a characterization of toric hyperk\"{a}hler varieties among symplectic singularities, 
Selecta Math. (2026) 32:4, 1-49 
\vspace{0.2cm}
 
[Lo] Losev, I.:  Classification of multiplicity free Hamiltonian actions of algebraic tori on Stein manifolds. (English summary)
J. Symplectic Geom. {\bf 7} (2009), no. 3, 295 - 310.
\vspace{0.1cm}

[Hu] Hubbard, A.: Classifying toric hyperk\"{a}hler varieties, in workshop ``Deformations and Birational Geometry of Algebraic Varieties'' , 
RIMS, Kyoto Univ. Jan. 19 -23, 2026

\vspace{0.5cm}

\begin{center}
Research Institute for Mathematical Sciences, Kyoto University, Oiwake-cho, Kyoto, Japan

E-mail address: namikawa@kurims.kyoto-u.ac.jp  
\end{center}

\end{document}